\documentclass[12pt,reqno]{amsart}

\usepackage{amsmath,amssymb,amsthm}
\usepackage{graphicx}
\usepackage[all]{xypic}
\usepackage[pdftex]{hyperref}
\input xypic
\newtheorem{problem}{Problem}[section]
\newtheorem{example}{Example}[section]
\newtheorem{remark}{Remark}[section]
\newtheorem{theorem}{Theorem}[section]
\newtheorem{corollary}{Corollary}[section]
\newtheorem{lemma}{Lemma}[section]                                                                                                                                                                                      \newtheorem{Proposition}{Proposition}[section]
\newtheorem{definition}{Definition}[section]

\numberwithin{equation}{section}
\def\nga{\widetilde}
\def\rem{\begin{remark}}
\def\erem{\end{remark}}
\def\ex{\begin{example}
  }
  \def\eex{\end{example}}
\def\thr{\begin{theorem}}
\def\ethr{\end{theorem}}
\def\pro{\begin{Proposition}}
\def\epro{\end{Proposition}}
\def\coro{\begin{corollary}}
\def\ecoro{\end{corollary}}
\def\df{\begin{definition}}
\def\edf{\end{definition}}
\def\lm{\begin{lemma}}
\def\elm{\end{lemma}}
\def\pf{\begin{proof}}
\def\epf{\end{proof}}
\def\problem{\begin{problem}}
\def\eproblem{\end{problem}}

\def\it{\begin{itemize}}
\def\hit{\end{itemize}}

\newcommand{\vm}{\forall}
\newcommand{\tru}{\setminus}

\newcommand{\sr}{\longrightarrow}

\newcommand{\can}{\sqrt}
\newcommand{\mtn}{\rightarrow}

\newcommand{\gt}{\overline}

\newcommand{\seq}[1]{\left<#1\right>}

\def\Coker{\mbox{Coker}}
\def\Ker{\mbox{Ker}}

\def\NRed{\mbox{NRed}}
\def\Sing{\mbox{Sing}}
\def\rad{\mbox{rad}}

\def\n{\Bbb N}

\def\c{\Bbb C}

\def\ohoa{\mathcal{O}}

\begin{document}
\title[Equinormalizable and equisingular deformations]
{Equinormalizable theory for  plane curve singularities with embedded  points and the theory of equisingularity}

\author{C\^{o}ng-Tr\`{i}nh L\^{e} }

\address{Department of Mathematics, Quy Nhon University\\
170 An Duong Vuong, Quy Nhon City, Vietnam}

\email{lecongtrinh@qnu.edu.vn}


\begin{abstract}
In this paper we give some criteria  for a family of generically reduced plane curve singularities to be equinormalizable. The first criterion is based on the $\delta$-invariant of a (non-reduced) curve singularity which is introduced by Br\"{u}cker-Greuel (\cite{BG}). The second criterion is based on the I-equisingularity of a $k$-parametric  family ($k\geq 1$) of generically reduced plane curve singularities, which is introduced by Nobile (\cite{No}) for one-parametric families. The equivalence of some kinds of equisingularities of a family of generically reduced plane curve singularities is also studied.

\end{abstract}
\maketitle


\section{Introduction}
The theory of equinormalizable deformations has been initiated by Teissier (\cite{Tei1}) for deformations of reduced curve singularities over $(\c,0)$. It is generalized to higher dimensional base spaces by Teissier himself and Raynaud in 1980 (\cite{Tei2}). Recently, it is developed by Chiang-Hsieh and Lipman (\cite{Ch-Li}) for projective deformations of reduced  complex spaces over normal base spaces, and it is studied by  Koll\'{a}r for projective deformation of generically reduced algebraic schemes over semi-normal base spaces (\cite{Ko}). The theory of equinormalizable deformations of not necessarily reduced curve singularity over $(\c,0)$ is studied by Br\"{u}cker and Greuel in 1990 (\cite{BG}). Some generalizations of the results of Br\"{u}cker and Greuel to deformations of (not necessarily reduced) curve singularities over normal base spaces are given by Greuel and the author  in the forthcoming  paper \cite{GL}. In this paper we study the equinormalizable deformations of not necessarily reduced plane curve singularities over smooth base spaces $(\c^k,0)$ ($k\geq 0$). We show in Theorem \ref{thr3.4.1} that  \emph{the induced morphism on the pure-dimensional part of the  total space is equinormalizable if and only if the given deformation is $\delta$-constant}. This result is generalized  in \cite{GL} to deformations of not necessarily reduced curve singularities over normal spaces, however the technique used to prove it in the context of plane curve singularities (in $\c^2$) is quite special, applying a consequence of the Hilbert-Burch Theorem (Lemma \ref{lm4.1}). This technique cannot be used for deformations of curve singularities which are not planar.

The theory of equisingularity for reduced plane curve singularities has been introduced  by Zariski (1970, \cite{Za}), Wahl (1974, \cite{Wa}). In \cite{No} Nobile  defined three kinds of equisingularities for a one-parametric family of generically reduced plane curve singularities: I-, T- and C-equisingularity. He showed that I-equisingularity is equivalent to T-equisingularity, C-equisingularity implies I-equisingularity, and gave a criterion for a family to be C-equisingular. In this paper we generalize the results of Nobile to $k$-parametric families ($k\geq 1$) of generically reduced plane curve singularities, and based on these equingularities we give a criterion for a deformation  to be equinormalizable (Theorem \ref{thr5.3}).

\section{$\delta$-invariant of (not necessarily reduced) curve singularities }
Following  Greuel and Br\"{u}cker (\cite{BG}, for curves), Greuel and the author (\cite{GL}, for arbitrary complex spaces), we recall in this section the definition of the $\delta$-invariant of a curve  which is not necessarily reduced, having an \emph{isolated singularity}\footnotemark. \footnotetext{A point $p$ in a curve  $C$  is said to be \emph{non-reduced} (resp. \emph{singular})  if the local ring at $p$, $\ohoa_{C,p}$ is not reduced (resp. not regular). The set of all non-reduced (resp. singular) points in $C$ is denoted by $\NRed(C)$ (resp. $\Sing(C)$),  and called the \emph{non-reduced locus} (\emph{singular locus})  of $C$. If $p\in \NRed(C)$ (resp. $\Sing(C)$)  is  isolated  then it is called an \emph{isolated non-reduced point} (resp.  \emph{isolated singular point})  of $C$. }

For a complex curve  $C$, we denote by $C^{red}$  its reduction and by  $i: C^{red} \mtn C$ the inclusion. For $\nu^{red}: \gt{C} \mtn C^{red}$ we mean the normalization of the reduced curve $C^{red}$, and we call the composition $\nu:=i \circ \nu^{red}: \gt{C} \mtn C $ the \emph{normalization of the curve $C$}. Then we have the induced map on the structure sheaves
$$\nu^\sharp:  \ohoa_C \mtn \nu_*\ohoa_{\gt{C}}. $$
We have $\Ker(\nu^\sharp)= Nil(\ohoa_C)$,  the sheaf of nilpotent elements of $\ohoa_X$, and $\Coker(\nu^\sharp)=\nu_*\ohoa_{\gt{C}}/\ohoa_{C^{red}}$. Since the map $\nu$ is finite, these sheaves   are  coherent $\ohoa_C$-modules, whose supports  are $\NRed(C)$ and  $\Sing(C)$, respectively. Thus, if $x\in C$ is an isolated non-normal point then $\Ker(\nu^\sharp)$ and $\Coker(\nu^\sharp)$ are finite dimensional  $\c$-vector spaces, and we have
$\Ker(\nu^\sharp) = H_{\{x\}}^0(\ohoa_C)$, the local cohomology.

\df \rm Let $C$ be a complex curve and $x\in C$ an isolated singular point. The number
$$\delta(C^{red},x):=\dim_\c (\nu^{red}_*\ohoa_{\gt{C}})_x/\ohoa_{C^{red},x} $$
is called the \emph{delta-invariant of $C^{red}$ at $x$},
$$\epsilon(C,x):=\dim_\c H_{\{x\}}^0(\ohoa_C) $$
is called  the \emph{epsilon-invariant of $C$ at $x$}, and
$$\delta(C,x):=\delta(C^{red},x) - \epsilon(C,x) $$
is called the \emph{delta-invariant of $C$ at $x$}.\\
If $C$ has only finitely many singular points then the number
$$\delta(C):=\sum_{x\in \Sing(C)} \delta(C,x) $$
is called the \emph{delta-invariant } of $C$.
\edf
It is easy to see that $\delta(C^{red},x)\geq 0$, and $\delta(C^{red},x) = 0$ if  and only if $x$ is an isolated point of $C$ or the germ $(C^{red},x)$ is smooth. Hence, if $x\in C$ is an isolated point of $C$ then $\delta(C,x) = -\dim_\c \ohoa_{C,x}$.

\ex \label{ex1} \rm We compute the $\delta$-invariant of the curve singularity $(X_0,0)\subseteq (\c^2,0)$ defined by the ideal
$$ I_0=\seq{x^2-y^3}\cap \seq{y}\cap\seq{x,y^5}\subseteq \c\{x,y\}. $$
The curve singularity $(X_0,0)$ is the union  of a  cusp $C$ and a straight line  $L$ with an embedded non-reduced point at the origin. We have 
$$ \delta(X_0^{red},0)=\delta(C,0)+\delta(L,0)+i_{(0,0)}(C,L) = 1 + 0 + 2 = 3,$$
where $i_{(0,0)}(C,L)$ denotes the intersection multiplicity of $C$ and $L$ at the origin.
Note that $H_{\{0\}}^0(\ohoa_{X_0}) = Nil(\ohoa_{X_0,0})$ is the kernel of the surjection
$$ p : \ohoa_{X_0} \cong \c\{x,y\}/I_0 \twoheadrightarrow
\ohoa_{X_0^{red}} \cong \c\{x,y\}/\rad(I_0).$$
Hence
$$\epsilon(X_0,0) = \dim_\c \rad(I_0)/I_0 = 1,$$
where $\rad(I_0) = \seq{x^2-y^3}\cap \seq{y}$.  Therefore
$$\delta(X_0) = \delta(X_0,0) = \delta(X_0^{red},0) - \epsilon(X_0,0) = 2.$$
\eex

\section{Simultaneous normalizations and equinormalizable morphisms}
Following Koll\'{a}r (\cite{Ko}) and Chiang-Hsieh-Lipman (\cite{Ch-Li}),  in \cite{GL} we  gave the following definition of a simultaneous normalization of a map between complex spaces, and we defined also equinormalizable maps between complex spaces (or germs).
\df \rm    Let $f: X\sr S$ be a morphism of complex spaces. A \emph{simultaneous normalization of $f$ } is a morphism
$n: \nga{X} \sr X$ such that
 \it
\item[(1)]  $n$ is finite,
\item[(2)] $\nga{f}:=f\circ n: \nga{X}\mtn S$ is \emph{normal}, i.e.,  for each $z\in \nga{X}$, $\nga{f}$ is flat at $z$ and the fiber $\nga{X}_{\nga{f}(z)}:=\nga{f}^{-1}(\nga{f}(z))$ is normal,
\item[(3)] the induced map $n_s: \nga{X}_s:=\nga{f}^{-1}(s) \sr X_s$ is \emph{bimeromorphic\footnotemark}\footnotetext{A map $f: X \mtn S$ is called \emph{bimeromorphic } if there exists a nowhere dense analytic subset $A$ of $S$ such that $f^{-1}(A)$ is nowhere dense in $X$ and the induced map  $X\tru f^{-1}(A) \mtn S\tru A$ is an isomorphism.} for each $s\in f(X)$.
\hit
The morphism $f$ is called \emph{equinormalizable} if  the normalization $\nu: \gt{X}\mtn X$ is a simultaneous normalization of $f$. It  is called \emph{ equinormalizable at $x\in X$} if the restriction of $f$ to some neighborhood of $x$ is equinormalizable.\\
If $f: (X,x) \sr (S,s)$ is a morphism  of germs, then a \emph{simultaneous normalization of $f$} is
a morphism $n$ from a multi-germ $(\nga{X}, n^{-1}(x))$ to $(X,x)$ such that some representative of $n$ is a simultaneous normalization of a representative of $f$. The germ $f$ is \emph{equinormalizable} if some representative of $f$ is equinormalizable.
 \edf
 \rem \rm Our definition of simultaneous normalizations of a map  does not require the flatness of the map. We do also not require the reducedness of the fibers (however, all non-empty fibers are generically reduced if the map admits a simultaneous normalization, see \cite{GL}). The total space of the map is also not required to be pure-dimensional. One may see \cite{Ch-Li}, \cite{Ko} and \cite{GL} for more disscussions about simultaneous normalizations and equinormalizability.
 \erem
 In this paper we consider flat deformations of plane curve singularities,  find a criterion  for such a deformation to be equinormalizable. First of all we give an example of a deformation of a plane  curve singularity that is equinormalizable.

 \ex  \label{ex2} \rm We consider again the curve singularity $(X_0,0)
\subseteq (\c^2,0)$ given in Example \ref{ex1}.  Let
$$f : (X,0) \sr (\c^2,0), (x,y,u,v) \mapsto (u,v), $$
 be the restriction to $(X,0) \subseteq (\c^4,0)$ of the projection
$\pi: (\c^4,0) \mtn (\c^2,0)$,  where  $(X,0)$ is defined by the ideal
$$I:=\seq{x^2-y^3+uy^2} \cap \seq{y-u}\cap \seq{x-v,y} \subseteq \c\{x,y,u,v\}.$$
The map $f$ is flat since $u,v$ is an $\ohoa_{X,0}$-regular sequence.
 Hence the map $f$  is a deformation of $(X_0,0)$ over $(\c^2,0)$.\\
It is easy to see that the  total space $(X,0)$ of the deformation $f$ is reduced with two $3$-dimensional irreducible components and one $2$-dimensional  irreducible component. The normalization of these  components are given by
\begin{align}
\nu_1 : (\c^3,0) & \sr (X,0)\nonumber\\
(T_1,T_2,T_3) & \mapsto (T_3^3+T_1T_3,T_3^2+T_1,T_1,T_2),\nonumber
\end{align}
\begin{align}
\nu_2 : (\c^3,0) & \sr (X,0)\nonumber\\
(T_1,T_2,T_3) & \mapsto (T_1, T_2,T_2,T_3),\nonumber
\end{align}
and 
\begin{align}
\nu_3 : (\c^2,0) & \sr (X,0)\nonumber\\
(T_1,T_2) & \mapsto (T_2, 0,T_1,T_2).\nonumber
\end{align}
For each $i=1,2,3$, the morphism $\bar{f}_i:=f\circ \nu_i $ is given by the last two components of $\nu_i$
and all of them are flat. The special fibers of $\bar{f}_1$ and $\bar{f}_2$ are straight lines in $\c^3$, while the special fiber of $\bar{f}_3$ is a (normal) point in $\c^2$.  It follows  that the given deformation is equinormalizable.
\eex
\section{Equinormalizability of deformations of curve singularities in the plane with embedded non-reduced points }
In this section we  give a criterion for a deformation of  an isolated curve singularity in the plane with embedded non-reduced points  over a smooth base space of dimension $k \geq 1$ to be equinormalizable. As a consequence of the Hilbert-Burch theorem (cf. \cite[Theorem 5]{Bur} or \cite[Satz 6.1]{BG}),  each ideal defining a curve singularity in $\c^n$ can be factorized as a product of an ideal defining a hypersurface singularity and an ideal defining a Cohen-Macaulay singularity of codimension $2$ in $\c^n$. For $n=2$,  the  hypersurface singularity is the pure $1$-dimensional part of the curve singularity  and   the Cohen-Macaulay singularity becomes a (non-reduced) point.

The following  result is one of the main ideas in the proof of  the numerical criterion for the equinormalizability.
\lm[cf. \cite{BG}, Prop. 6.3] \label{lm4.1} Let $(X_0,0) \subseteq (\c^n,0)$ be a
curve singularity defined by the ideal $\seq{g_1,\cdots,g_m} =
\seq{g}.\seq{p_1,\cdots,p_m}.$ Let $f : (X,0) \sr (\c^k,0)$ be a
deformation of $(X_0,0)$, where $(X,0)$ is given by the ideal $I(X,0) =
\seq{g_i+ \sum_{j=1}^k t_j \bar{g}_{ij}}_{1\leq i \leq m} \subseteq
\ohoa_{\c^n,0}\{\textbf{t}\}$, $\textbf{t}:= (t_1,\cdots,t_k)$. Then
there exist functions $\bar{g}_{j}$ and $\bar{p}_{ij}$ in
$\ohoa_{\c^n,0}\{\textbf{t}\}$ such that
$$ \seq{gp_i+ \sum_{j=1}^k t_j \bar{g}_{ij}}_{1\leq i \leq m} =
\seq{g + \sum_{j=1}^k t_j \bar{g}_{j}}\cdot \seq{p_i+ \sum_{j=1}^k
t_j \bar{p}_{ij}}_{1\leq i \leq m},$$ where $f_G : (G,0) \sr
(\c^k,0)$ with $I(G,0) = \seq{g + \sum_{j=1}^k t_j \bar{g}_{j}}$ is a
deformation of $(G_0,0)$ defined by $\seq{g}\subseteq \ohoa_{\c^n,0}$,
and $f_P : (P,0) \sr (\c^k,0)$ with $I(P,0) = \seq{p_i+ \sum_{j=1}^k
t_j \bar{p}_{ij}}_{1\leq i \leq m}$ is a deformation of $(P_0,0)$
defined by $\seq{p_1,\cdots,p_m}\subseteq \ohoa_{\c^n,0}$.
\elm
Let  $f: X\mtn S$ be a morphism of complex spaces whose fibers have only finitely many non-normal points. It is called \emph{(locally) delta-constant} if the function $s \mapsto \delta(X_s)$ is (locally) constant on $S$. A morphism of gemrs is \emph{$\delta$-constant } if some of its representatives is $\delta$-constant.

 The following theorem is the first main result of this paper.
 \thr \label{thr3.4.1}
Let $(X_0,0) \subseteq (\c^2,0)$ be an isolated (not necessarily
reduced) curve singularity. Let $f : (X,0) \sr (\c^k,0), k \geq 1$,
 be a deformation of $(X_0,0)$.  Denote by $(X^u,0)$ the \emph{unmixed subgerm\footnotemark} \footnotetext{Let $R$ be a ring and  $I\subseteq R $ be an ideal of dimension $m$.  Assume that $I$ has an  irredundant primary decomposition $I=\bigcap_{i=1}^r Q_i$. For an integer  $0\leq k \leq m$, we define the \emph{pure $k$-dimensional part} $I^{(k)}$ of the ideal $I$ to be the intersection of all $Q_i$  with $\dim Q_i = k$. The ideal  $I^{(k)}$ is well-defined for each $0\leq k \leq m$ because this part of the primary decomposition is uniquely determined. We define the \emph{unmixed part} $I^u$ of the ideal $I$ to be the pure $m$-dimensional part of the radical $\can{I}$. If the germ $(X,0)\subseteq (\c^n,0)$ is defined by an ideal $I\subseteq \ohoa_{\c^n,0}$, the \emph{unmixed subgerm } $(X^u,0)$ of $(X,0)$ is the one defined by the unmixed part $I^u$ of $I$.  } of $(X,0)$.  Then the following holds:
  \it
  \item[(1)] The restriction $f^u: (X^u ,0) \sr (\c^k,0)$\index{$f^u$} \index{$X^u$} is flat.
  \item[(2)] $f$ is $\delta$-constant  if and only if  $f^u: (X^u,0)\sr (\c^k,0) $ is equinormalizable.
\hit
 \ethr
\pf
 (1) As a consequence of the Hilbert-Burch theorem, the germ  $(X_0,0)$ is a union
of a hypersurface $(X_0^{u},0)\subseteq (\c^2,0)$   and an embedded  (non-reduced) point. More precisely, if $(X_0,0)$ is
defined by the ideal $I_0$ then it can be factorized as 
 $$I_0 = \seq{g}\cdot J_0,$$
 where $\seq{g} \subseteq \ohoa_{\c^2,0}$ defines
$(X_0^{u},0)\subseteq (\c^2,0)$ and the ideal  $J_0\subseteq
\ohoa_{\c^2,0}$ defines such a (non-reduced) point. By Lemma \ref{lm4.1}
we can write the ideal $I$ defining $(X,0)\subseteq (\c^2\times
\c^k,0)$ as  $ I  = \seq{G}\cdot J$, where $\seq{G}\subseteq \ohoa_{\c^2\times \c^k,0}$ defines a
deformation $(H,0) \subseteq (\c^2\times \c^k,0)$ of the hypersurface
$(X_0^{u},0)$ given by $\seq{g}$ and $J\subseteq \ohoa_{\c^2\times
\c^k,0}$ defines a deformation $(P,0)\subseteq (\c^2\times \c^k,0)$ of
the (non-reduced) point given by $J_0$. \\
 Note that $(H,0)$ is reduced and pure $(k+1)$-dimensional, because it is the
total space of a deformation of the reduced  and pure 1-dimensional  singularity
$(X_0^{u},0)\subseteq (\c^2,0)$ over $(\c^k,0)$ which is  reduced and pure $k$-dimensional (cf. \cite[Theorem I. 1.85]{GLS} and \cite[Theorem I. 1.101]{GLS}).  Hence $(H,0)
\equiv (X^u,0)$. Therefore the restriction map $f^{u}:
(X^u,0) \equiv (H,0) \sr (\c^k,0)$ is flat and it is actually a deformation of $(X_0^u,0)$ over $(\c^k,0)$.\\
(2) \textbf{First we prove the "only if" part}.  Let $f : X \sr S$ be a sufficiently small representative of the given deformation such that $f^u: X^u \sr S$  is equinormalizable. Since $(X_0,0)$ has isolated singularities, it follows from the \emph{generic principle} (cf. \cite[Theorem 2.2]{BF}) that there exists an open dense subset $U \subseteq S$ such that $(X^{u})_s:=(f^u)^{-1}(s)$ are reduced  for all $s \in U$. \\
We first show  that \emph{$f$ is $\delta$-constant on $U$}, i.e., $\delta(X_s) = \delta(X_0)$ for any $s \in U$.  In fact, for any $s\in U$,  there exist an irreducible reduced cure singularity $C \subseteq S$ passing through $s$. Let $\alpha: T \sr C \subseteq S$ be the normalization of this curve singularity such that $\alpha(T\tru \{0\}) \subseteq U, \alpha(0)=s$, where $T\subseteq \c$ is a small disc with center at $0$. Denote
$X_T\index{$X_T$}:=X\times_S T,~ X^u_T\index{$X^u_T$}:=X^u\times_S T, ~ \gt{X}_T\index{$\gt{X}_T$}:= \gt{X}\times_S T,$
 where $\nu^u: \gt{X} \mtn X^u$ is the normalization of $X^u$. Thus  we have the following Cartesian diagram:
 $$\xymatrix@C=12pt@R=8pt@M=6pt{
&&&\ar @{} [dr] |{\Box} \gt{X}_T \ar[r] \ar[d]_{\nu^u_T} \ar@/_3pc/[ddd]_{\bar{f}_T} & \gt{X} \ar[d]^{\nu^u} \ar@/^3pc/[ddd]^{\bar{f}}\\
\ar @{} [dr] |{(\Delta)}&&&\ar @{} [dr] |{\Box} X^u_T \ar[r] \ar @{^{(}->}[d]^{i_T} \ar@/_2pc/[dd]^{f^u_T} & X^{u} \ar @{^{(}->}[d]_i \ar@/^2pc/[dd]_{f^u\index{$f^u$}}\\
&&&\ar @{} [dr] |{\Box} X_T\ar[d]^{f_T} \ar[r]& X \ar[d]_f\\
&&&T \ar[r] & S}$$
For any $t\in T, s = \alpha(t) \in S$, we have
$$\ohoa_{(X_T)_t}:= \ohoa_{f_T^{-1}(t)} \cong \ohoa_{X_s},$$
 $$\ohoa_{(X^u_T)_t}:= \ohoa_{{f^u_T}^{-1}(t)} \cong \ohoa_{(X^u)_s},~\ohoa_{(\gt{X}_T)_t}:= \ohoa_{\bar{f}_T^{-1}(t)} \cong \ohoa_{\gt{X}_s}.$$
In the diagram $(\Delta)$, $f$ is flat by hypothesis, $\bar{f}$ is flat since $f^u$ is equinormalizable and  $f^u$ is flat by (1). It follows from the preservation of flatness under base change (cf. \cite[Prop. I. 1.87]{GLS}) that the induced morphisms $f_T, \bar{f}_T$ and $f^u_T$ are also flat over $T$. Moreover, for any $t\in T\tru \{0\}$ we have $s:=\alpha(t)\in U$, hence $(X^u_T)_t\cong (X^u)_s$ is reduced. It follows from \cite[Prop. 3.1.1 (3)]{BG} that $X^u_T$ is reduced. On the other hand, since $X^u$ and $\c^k$ are pure dimensional, the fiber of $f^u$ is pure dimensional. It implies that the fiber of $f^u_T: X^u_T\mtn T$ is also pure dimensional. Hence $X^u_T$ is pure dimensional. Thus $X^u_T$ is the unmixed space $(X_T)^u$ of $X_T$. \\
 Moreover, since $(\gt{X}_T)_t\cong \gt{X}_s$ for any $t \in T, s = \alpha(t),$  it implies that the special fiber $(\gt{X}_T)_0$ of $\bar{f}_T$ is normal. Then $\bar{f}_T$ is regular by the \emph{regularity criterion for morphisms} (cf. \cite[Theorem I. 1.117]{GLS}).  It implies that $\gt{X}_T \cong (\gt{X}_T)_0 \times T$ which is smooth, hence normal, and it is the normalization of $X^u_T \equiv (X_T)^u$. Consider the morphism $f^u_T: X^u_T \sr T$ with $X^u_T$ reduced and pure 2-dimensional, hence $X^u_T$ is unmixed. Since $(\gt{X^u_T})_0 \cong (\gt{X}_T)_0 \cong \gt{X}_0$ which is normal, it implies that the morphism $f^u_T$ is equinormalizable. It follows from \cite[Korollar 2.3.5]{BG} that $f_T: X_T \sr T$ is $\delta$-constant, hence $f: X \sr S$ is $\delta$-constant on $U$. \\
Let us now take  $s_0 \in S\tru U$. Since $U$ is dense in $S$, $s_0 \in S$,
there exists always a  point $s_1 \in U$ which is closed to $s_0$.  It follows from the
semi-continuity of the $\delta$-function (Lemma \ref{lmDeltaSemi}) that
$$ \delta(X_0) \geq \delta(X_{s_0}) \geq \delta(X_{s_1}).$$
Moreover,  $\delta(X_0) = \delta(X_{s_1})$ as above. It implies that $\delta(X_{s_0})=\delta(X_0)$.  Hence $f : X \sr S$ is $\delta$-constant. \\
\textbf{Now we show the "if" part}. Let $f : X \sr S$ be a sufficiently small representative of the given deformation such that it is $\delta$-constant.   For each $s \in S$, denote $X_s:= f^{-1}(s)$  and for each $x_s
\in \Sing(X_s)$, consider the family
$$\xymatrix@C=12pt@R=12pt@M=6pt{
  (X_s,x_s)\ar[d] \ar@{^{(}->}[r] & (X,x_s)\ar[d] \\
\{s\}\ar@{^{(}->}[r]& (S,s) }$$
Each germ  $(X_s,x_s) \subseteq (\c^2,0)$ is defined by an ideal of the form
$\seq{g_s}\cdot J_s \subseteq \ohoa_{\c^2\times\{s\},x_s}.$ We have
$$\delta(X_s,x_s) = \delta(X_s^{u},x_s) - \dim_\c
\ohoa_{\c^2\times\{s\},x_s}/J_s,$$ and
\begin{equation} \label{equ1}
\begin{array}{l}\delta(X_s) := \sum_{x_s\in
\Sing(X_s)}\delta(X_s,x_s)\\
  = \sum_{x_s\in \Sing(X_s)} \delta(X_s^{u},x_s) - \sum_{x_s\in \Sing(X_s)} \dim_\c
\ohoa_{\c^2\times\{s\},x_s}/J_s.
\end{array}
\end{equation}
Since $(X_0,0)$ is isolated, it follows from the \emph{local finiteness theorem} (cf. \cite[Theorem
I.1.66]{GLS}) that  the restriction  $f: \Sing(f) \sr S$ is finite and $\Sing(X_0)=\Sing(f)\cap X_0 = \{0\}$. Then  $f : P \sr S$ is also finite. Obviously, it is flat. Therefore it follows from the  \emph{semi-continuity of fibre
functions} (cf. \cite[Theorem I.1.81]{GLS}) that for all $s \in S$, we have
\begin{equation}\label{equ2}
 \dim_\c\ohoa_{\c^2\times\{0\},0}/J_0 = \sum_{x_s\in \Sing(X_s)} \dim_\c
\ohoa_{\c^2\times\{s\},x_s}/J_s.
\end{equation}
Moreover, $f$ is $\delta$-constant  by assumption, that is,
\begin{equation}\label{equ3}
 \delta(X_0) = \delta(X_s) \mbox{ for all } s \in S.
\end{equation}
It follows from (\ref{equ1}), (\ref{equ2}) and (\ref{equ3}) that
$$\delta(X_0^{u}) = \delta(X_0^{u},0) =   \sum_{x_s\in \Sing(X_s)}
\delta(X_s^{u},x_s) = \delta(X_s^{u})
 \mbox{ for all } s \in S,$$
i.e., $\delta(X_s^{u})$ is constant. Therefore we have a
\emph{delta}-constant family of reduced curve singularities $f^{u}:
X^{u} \equiv H  \sr S$. Then it is equinormalizable  by the criterion of Teissier, Raynaud, Chiang-Hsieh and Lipman (cf. \cite[Theorem 5.6]{Ch-Li}).
 \epf
 In the proof of the theorem above we used the following semi-continuity of the delta-function.
 \lm \label{lmDeltaSemi} Let $f: (X,0) \sr (S,0)$ be a deformation of an isolated (not necessarily reduced) curve singularity $(X_0,0)\subseteq (\c^n,0)$. Then the $\delta$-function, $s \mapsto \delta(X_s)$, is \emph{upper semi-continuous}  in the following sense: there exists a representative $f : X \sr S$ of the given deformation such that  $\delta(X_s) \leq
\delta(X_0)$  for all $s \in S$.
\elm

\pf Let $f : X \sr S$ be a sufficiently small representative of the
given deformation. For any $s \in S$, there exists an irreducible  reduced curve
singularity $C\subseteq S $ passing through $0$ and  $s$ and let $\alpha: T \sr
C \subseteq S $ be the normalization of this curve singularity, where
$T \subseteq \c$ is a small disc.  Denote
$ X_T\index{$X_T$}:=X\times_S T$. Let  $f_T: X_T \mtn T$  be the morphism induced by $f$. The morphism  $f_T$ is flat by the preservation of flatness under base change (cf.  \cite[Prop. I.1.85]{GLS}). Moreover, for any $t \in T, s:=\alpha(t)\in S$, we have
$ \ohoa_{(X_T)_t}:= \ohoa_{f_T^{-1}(t)}\cong \ohoa_{f^{-1}(s)}=:\ohoa_{X_s}$.
It follows from \cite[Satz 3.1.2 (iii)]{BG} that, for $t\in T$
such that $\alpha(t) = s$, we have
$$ \delta(X_0) - \delta(X_s) = \delta((X_T)_0) - \delta((X_T)_t) \geq 0.$$

\epf

\ex \label{ex3} \rm Let us consider again the deformation $f: (X,0)\mtn (\c^2,0)$ of the plane curve singularity $(X_0,0)$ given in Example \ref{ex2}. As we have shown there, this deformation is equinormalizable. The $\delta$-invariant of the special fiber $(X_0,0)$ is equal to $2$ (Example \ref{ex1}). Moreover, for each $u\not = 0, v\not = 0$ close to $0$, the reduced fiber $X_{uv}=f^{-1}(u,v)$ consists of a cubic curve $C$, a straight line $L$ and a reduced point $(v,0)$. We can compute $\delta(X_{uv})= 2$. \\
Furthermore, the fibers $X_{0u}$ and $X_{v0}$ are non-reduced and their delta-invariants are also $2$. Hence the given deformation is $\delta$-constant. 
\eex

\rem \rm  With the above notations, if the total space $(X,0)$ of the deformation $f : (X,0) \sr (\c^k,0)$ of the plane curve singularity $(X_0,0)$  is reduced and  pure $(k+1)$-dimensional then $(X_0,0)$ is necessarily  reduced. In fact, since   $(X,0)$ is reduced,  the ideal $I$ defining $(X,0)$ is radical, i.e., $I=\can{I}$. Moreover, since  $(H,0)$ is
pure $(k+1)$-dimensional and $(P,0)$ is pure $k$-dimensional, it follows that $(P,0)\subseteq (H,0)$. Then   $(X,0) = (H,0) \cup (P,0) = (H,0)$, i.e., $V(I) = V(G)$. It follows from Hilbert-R\"{u}ckert's
Nullstellensatz (cf. \cite[Theorem I.1.72]{GLS}) that
$$\seq{G}\cdot J = I = \can{I} = I(V(I)) = I(V(G)) = \can{\seq{G}} =
\seq{G}.$$
 Hence $G \in \seq{G}\cdot J$. Then there exists $h \in J$
such that $G = Gh$, or $G(1-h)= 0$. Since $G$ is a non-zerodivisor of
$\ohoa_{\c^2\times \c^k,0}$ we get $h = 1$. Hence $ 1 \in J$ and we
have $J = \seq{1}$. This implies $I = \seq{G}$ and hence $I_0 =
\seq{g}$. It means that  $(X_0,0)$ is reduced. \\
Thus, if the plane curve singularity $(X_0,0)$ is not reduced then the total spaces of deformations  over smooth base spaces are either not reduced or not pure dimensional.
However this fact is not true for deformations of non-plane curve singularities. We show in the following example that there exists a   deformation of a non-reduced  curve singularity in a $4$-dimensional complex space which has a reduced and pure dimensional total space.
\erem

\ex \label{ex4} \rm Let us consider the curve singularity  $ (X_0,0)\subseteq (\c^4,0)$ defined by the ideal 
$$I_0:= \seq{x^2 - y^3,z,w} \cap \seq{x,y,w} \cap \seq{x,y,z,w^2} \subseteq \c\{x,y,z,w\},$$
 which  was considered by  Steenbrink (\cite{St}). The curve singularity $(X_0,0)$ is a union of a cusp
$C$ in the plane $z=w=0,$ a straight line $L = \{x = y = w = 0\}$ and
an embedded non-reduced point $O = (0,0,0,0)$. Now we consider the restriction $f: (X,0)\mtn (\c^2,0)$ of the projection $\pi:(\c^6,0)\mtn (\c^2,0), ~ (x,y,z,w,u,v)\mapsto (u,v),$ to the complex space $(X,0)$ which is defined by the ideal 
$$I=\seq{x^2-y^3+uy^2,z,w} \cap \seq{x,y,w-v}\subseteq \c\{x,y,z,w,u,v\}.$$
It is easy to show that the total space $(X,0)$ is reduced and pure $3$-dimensional.  Moreover, by a similar way to Example \ref{ex2} and Example \ref{ex3} we can show that this deformation is equinormalizable and delta-constant (with $\delta=1$).
\eex

\section{The theory of equisingularity }
In this section we study the theory of equisingularity for plane curve singularities with embedded points which is introduced in \cite{No} (for one-parametric family), where the author formulated and proved the equivalence between I-equisingularity and T-equisingularity, also the relation between I-equisingularity and C-equisingularity.

\df \rm A \emph{$k$-parametric family} ($k\in \n^*$) of generically reduced plane curve singularities is a diagram
$$
\xymatrix@C=8pt@R=8pt@M=4pt{
  (X,0)\ar@{^{(}->}[rr]^i \ar[dr]_f& & (\c^{k+2},0)\ar[dl]^\pi \\\
  & (\c^k,0),& }
$$
where $\pi$ is smooth and surjective, $f$ is flat, $X_t:=p^{-1}(t) $ is  a generically plane curve singularity for each $t\in \c^k$ closed to $0$. We denote this family shortly by $(X,f,\pi)$. As we have seen in the previous sections, this family is a \emph{deformation } of the generically reduced plane curve singularity $(X_0,0)$, which is the \emph{special fiber } of the deformation.  Throughout this section we restrict our attention to the families whose restriction of $f$ on its singular locus\footnotemark \footnotetext{The \emph{singular locus} of a flat map $f: X \mtn S$ is the set of all points $x\in X$ such that the fiber $X_{f(x)}:=f^{-1}(f(x))$ is singular. If $S$ is regular and $f$ is flat then $\Sing(X)\subseteq \Sing(f)$ (\cite[Theorem I.1.117]{GLS}).}  $ \Sing(f)$ is finite.
\edf
Assume that the germ $(X,0)$ is defined by the ideal
$$I:=I(X,0)\subseteq \c\{x,y,u_1,\cdots,u_k\},$$
 where $x,y$ (resp.  $u_1,\cdots,u_k$) are local coordinates in $\c^2$ (resp. $\c^k$). Then a $k$-parametric family of plane curve singularities induces a $k$-parametric flat family of plane ideals $(I,\pi)$ (compare to \cite[Definition 3.3]{No}).

We  associate to  each ideal $J \subseteq \c\{x,y\}$ a \emph{weighted  directed tree} $\tau(J)$ as defined in \cite[Section 1.2]{No}. The induced family of plane ideals  $(I,\pi)$ mentioned  above is called \emph{equisingular } if $\tau(I(t)) \thickapprox \tau(I(0))$ for all $t\in \c^k$ closed to $0$, where $I(t):=I\ohoa_{X_t,0}\subseteq\c\{x,y\}$.

\df \rm A $k$-parametric family of generically reduced plane curve singularities $(X,f,\pi)$  is said to be \emph{I-equisingular } if the induced family of plane ideals $(I:=I(X),\pi)$ is equisingular in the sense of weighted directed tree.
\edf
We also associate to each generically reduced plane curve $C\subseteq \c^2$ a \emph{bi-weighted directed tree } $T_2(C,\gamma)$ as defined in \cite[Section 2.4]{No}, where $\gamma$ is an ordering of the branches of the reduction $C^{red}$ of $C$.

\df \rm A $k$-parametric family of generically reduced plane curve singularities is said to be \emph{T-equisingular } if for any pair of points $t, t'$ in the same connected component of $\c^k$, closed to $0$, we can choose suitable orderings $\gamma, \gamma'$ on the reduction $X_t^{red}$ and $X_{t'}^{red}$ of the fibers $X_t, X_{t'}$ respectively such that the corresponding bi-weighted directed trees $T_2(X_t,\gamma_t)$ and $T_2(X_{t'},\gamma_{t'})$ are isomorphic.
\edf

\df \label{defC-equi} \rm A $k$-parametric family of generically reduced plane curve singularities is said to be \emph{C-equisingular } if it is I-equisingular and  if we denote by   $\pi_i: Z_i \mtn Z_{i-1}$ ($i\geq 1$) the blowing up of $Z_{i-1}$ with the center $\Sing(X_i)$, where  $X_i$ denotes the proper transform of $X$ under $\pi^{(i)}:=\pi\circ \pi_1 \circ \cdots \circ \pi_i: Z_i \mtn Z_0=Z$, then the induced morphism $(X_i,p_i) \mtn (\c^k,0)$ is flat, $p_i\in (\pi^{(i)})^{-1}(0)$.
\edf

The following theorem  gives an  equivalence of three kinds of equisingularity mentioned above for a $k$-parametric family of plane curve singularities. A similar result for one-parametric families is given by Nobile (\cite[Theorem 5.5 and Prop. 5.8]{No}).

\thr \label{thr5.1} Let $(X,f,\pi)$ be a  $k$-parametric family of  generically reduced plane curve singularities. Then
\it
\item[(i)] I-equisingularity is equivalent to T-equisingularity.
\item[(ii)] C-equisingularity implies I-equisingularity. Conversely, if the given family  is I-equisingular and all the fibers $X_t, t\in \c^k$ close to $0$, are smooth, then the family is C-equisingular.
\hit
\ethr

\pf It suffices to prove the theorem for a sufficiently small representative
$$
\xymatrix@C=8pt@R=8pt@M=4pt{
  X\ar@{^{(}->}[rr]^i \ar[dr]_f& & Z\ar[dl]^\pi \\\
  & S,& }
$$
of  the given family, where $X$, $Z$ and $S$ are sufficiently small neighborhoods of $0$.\\
 For each $s\in S$, there exists an irreducible reduced curve singularity $C\subseteq S$ passing through $s$. Let $\alpha: T \mtn S$ be the normalization of this reduced curve singularity, where $T\subseteq \c$ is a small disc with center at $0\in \c$. Denote
 by $X_T:=X\times_{S} T$, the Castesian product of $X$ and $T$ over $S$, and by $f_T: X_T \mtn T$ the induced morphism of $f$. By the preservation of flatness under base change (cf. \cite[Prop. I.1.87]{GLS}), $f_T$ is flat. Moreover, for each $t\in T$,  $s=\alpha(t)$, we have
 $$\ohoa_{(X_T)_t}= \ohoa_{X_T} \otimes_{\ohoa_{T,t}} \c = (\ohoa_X \otimes_{\ohoa_S,s} \ohoa_{T,t})  \otimes_{\ohoa_{T,t}} \c) \cong \ohoa_X \otimes_{\ohoa_{S,s}} \c \cong \ohoa_{X_s}. $$
Hence the fiber of $f$ over each $s\in S$ is isomorphic to the fiber of $f_T$ over $t\in T, \alpha(t)=s$.\\
 (i) The family $(X,f,\pi)$ is I-equisingular if and only if $\tau(I(s)) \thickapprox \tau(I(0))$ for all $s\in S$. Hence the trees    $\tau(I_T(t))$ and $ \tau(I_T(0))$ are isomorphic for all $t\in T$, where $I_T:=I_T(X_T)\subseteq \c\{x,y,t\}$ denotes the ideal defining $X_T$. Equivalently, the induced family $(X_T,f_T,\pi_T)$ is I-equisingular (by a theorem of Risler, cf.\cite[Theorem 3.7]{No}). Thus the family $(X_T,f_T,\pi_T)$  is T-equisingular. It means that the tree $T_2\big((X_T)_t,\gamma_t\big) $ is isomorphic to the tree $T_2\big((X_T)_{t'},\gamma_{t'}\big)$ for each pair $t,t'$ in the same connected component of $T$. Hence, for each pair $s,s'$ in the same connected component of $S$, $\alpha(t)=s, \alpha(t')=s'$,  the trees $T_2\big(X_s,\gamma_s\big) $ and $T_2\big(X_{s'},\gamma_{s'}\big) $ are isomorphic. This is equivalent to the T-equisingularity of the family $(X,f,\pi)$.\\
 (ii) It is clear that C-equisingularity implies I-equisingularity. Now we assume that the family is I-equisingular and all fibers $X_s, s\in S, $ are smooth. Since $f: X \mtn S$ is flat, $S$ is smooth, it follows that $X$ is smooth (cf. \cite[Theorem I.1.117]{GLS}). Therefore, for each $i \in \n^*$, the proper transform $X_i$ of $X$ under $\pi^{(i)}:=\pi\circ \pi_1\circ \cdots \circ \pi_i: Z_i \mtn Z_0=Z$ is smooth, of pure $(k+1)$-dimensional. Moreover,  each fiber $(X_i)_s:=(\pi^{(i)}\mid_{X_i})^{-1}(s)$ is a proper transform of the smooth fiber $X_s$, hence $(X_i)_s$ is smooth of pure $1$-dimensional. Thus we have the \emph{dimension formula }
 $$\dim (X_i,x) = \dim ((X_i)_s,x) + \dim(S,s), \pi^{(i)}(x)=s. $$
It follows that $\pi^{(i)}\mid_{X_i}$ is open (cf. \cite[Section 3.10, Theorem, p.145]{Fi}). Moreover, $X_i$ is smooth, hence Cohen-Macaulay. It follows that $\pi^{(i)}\mid_{X_i}$ is flat (\cite[Section 3.20, Proposition, p.158]{Fi}). Hence the given family is C-equisingular.
\epf
In the following we may induce an I-equisingular (hence T-equisingular) family of reduced plane curve singularities from a given I-equisingular family of generically reduced plane curve singularities.

\thr \label{thr5.2} Let $(X,f,\pi)$ be a $k$-parametric family of generically reduced plane curve singularities which is I-equisingular. Then all the fibers of the restriction $f^u: (X^u,0) \mtn (\c^k,0)$ are reduced and  the induced family
$$
\xymatrix@C=8pt@R=8pt@M=4pt{
  (X^u,0)\ar@{^{(}->}[rr]^i \ar[dr]_{f^u} & & (\c^{k+2},0)\ar[dl]^\pi \\\
  & (\c^k,0),& }
$$
is I-equisingular.
\ethr

\pf Let $f: X \mtn S$ be a sufficiently small representative of the germ $f: (X,0) \mtn (\c^k,0)$. Suppose $f$ is a deformation of the generically reduced plane curve singularity $X_0\subseteq \c^2$ which is defined by the ideal $I_0$. As a consequence of the Hilbert-Burch theorem, $I_0=\seq{g}\cdot J_0$, where $\seq{g}$ defines the unmixed part $X_0^u$ of $X_0$ and $J_0$ defines the embedded non-reduced point $0\in X_0$. It follows from Lemma \ref{lm4.1} and the proof of Theorem \ref{thr3.4.1} that the restriction $f^u: X^u \mtn S$ is a deformation of the reduction $X_0^u$, which is reduced and pure 1-dimensional. Hence the special fiber and the nearby fibers  of $f^u$ are all reduced.\\
Now we  consider the induced family
$$
\xymatrix@C=8pt@R=8pt@M=4pt{
  (X^u,0)\ar@{^{(}->}[rr]^i \ar[dr]_{f^u} & & (\c^{k+2},0)\ar[dl]^\pi \\\
  & (\c^k,0).& }
$$
Let $f^u : X^u \mtn S$ be a sufficiently small representative of the restriction map-germ $f^u$ in the family. By the same notation as in the proof of Theorem \ref{thr5.1} we have the following diagram
$$
\xymatrix@C=10pt@R=10pt@M=6pt{
  (X^u)_T\ar[r] \ar[d]_{(f^u)_T} &  X^u\ar[d]^{f^u} \\\
  T \ar[r]^{\alpha} & S,}
$$
where the induced map $f^u_T: (X^u)_T \mtn T$ is flat by the preservation of flatness under base change.
We have already showed in the proof of Theorem \ref{thr3.4.1} that $(X^u)_T=(X_T)^u$. Since the induced family $(X_T,f_T,\pi_T)$ is I-equisingular, it follows from \cite[Theorem 5.11]{No} that the family $\big((X_T)^u=(X^u)_T,f_T,\pi_T\big)$ is I-equisingular with reduced fibers. On the other hand, the fibers of $f^u_T: (X^u)_T \mtn T$ and $f^u: X^u \mtn S$ are isomorphic. Therefore the family $(X^u,f^u,\pi)$ is I-equisingular.
\epf

As a consequence we have the following equivalence of I-equisingularity and the delta-constancy of a family of generically reduced plane curve singularities. A family $(X,f,\pi)$ is said to be \emph{$\delta$-constant} if the morphism $f: (X,0) \mtn (\c^k,0)$ is $\delta$-constant.

\thr \label{thr5.3}  Let $(X,f,\pi)$ be a $k$-parametric family of generically reduced plane curve singularities. Then it is $\delta$-constant if and only if  the family  $(X^u,f^u,\pi)$  is I-equisingular.
\ethr

\pf The  I-equisingularity of the family  $(X^u,f^u,\pi)$ implies that  the induced family $\big((X^u)_T, f^u_T,\pi_T\big)$ is also I-equisingular with reduced fibers. It follows that the deformation $f^u_T: (X^u)_T \mtn T$ with reduced fibers is equinormalizable (cf. \cite[Theorem 5.3.1]{Tei2}). By a result of Teissier (cf. \cite[Corollary 1, p.609]{Tei1} ),
$$\delta\big(((X^u)_T)_t\big) = \delta\big(((X^u)_T)_0\big), \vm t \in T \mbox{ close to } 0.  $$
Since the fibers of $f^u: X^u \mtn S$ and $f^u_T: (X^u)_T \mtn T$ are isomorphic, it implies that the deformation $f^u : X^u \mtn S$ is $\delta$-constant. Hence it is equinormalizable (cf. \cite[Theorem 5.6]{Ch-Li}). It follows again from Theorem \ref{thr3.4.1} that $f: X \mtn S$ is $\delta$-constant. The converse of the argument given above is also satisfied. This proves the theorem.
\epf
\textbf{Acknowledgement:} {The author is indebted to Professor Gert-Martin Greuel for suggesting the topic of this paper and his constant help. He also thanks the referees for their useful comments and suggestions. }



\begin{thebibliography}{}
\bibitem  [BF]{BF} Bingener J. and  Flenner H., \emph{On the fibers of analytic mappings.} in Complex Analysis and Geometry, V. Ancona and A. Silva Eds.,  Plenum Press, 1993, pp. 45-102.
\bibitem  [BG]   {BG} Br\"{u}cker  C. and Greuel G.-M.,
\emph{Deformationen isolierter Kurven Singularit\"{a}ten mit
eigebetteten Komponenten.} (German) [Deformations of isolated curve
singulartities with embedded components] Manuscripta Math.
\textbf{70}, no. 1 (1990), 93-114.
\bibitem  [BuG]   {BuG} Buchweitz  R.O. and Greuel G.-M.,
\emph{The Milnor number and deformations of complex curve
singularities.} Inv. Math. \textbf{58} (1980), 241-281.
\bibitem  [Bur]   {Bur} Burch L.,
\emph{On ideals of finite homological dimension in local rings.}
Proc. Cam. Phil. Soc. \textbf{64} (1968), 941-948.
\bibitem  [Ch-Li]   {Ch-Li} Chiang-Hsieh H. J. and Lipman J., \emph{A numerical criterion for simultaneous
normalization.} Duke Math. J. \textbf{133}, no. 2 (2006), 347-390.
\bibitem  [Fi]   {Fi} Fischer G., \emph{Complex analytic geometry.} Lecture Notes in Math. \textbf{538}, Springer, 1976.
\bibitem  [GL]   {GL} Greuel G.-M. and  Le C.-T.,
\emph{On equinormalizable deformations of isolated singularities.}
In preparation, 2010.
\bibitem  [GLS]   {GLS} Greuel G.-M., Lossen C. and Shustin E.,
\emph{Introduction to Singularities and Deformations.} Springer, 2007.
\bibitem  [Ko]   {Ko} Koll\'{a}r  J., \emph{Simultaneous normalization and algebra husks.} Preprint, 2009. http://arxiv.org/abs/0910.1076
\bibitem  [No]   {No} Nobile A., \emph{Equsingularity theory for plane curves with embedded points.} Pacific J. Math. Vol. \textbf{170}, no. 2 (1995), 543-566.
\bibitem  [St]   {St} Steenbrink J., \emph{On mixed Hodge-structures.} Manuskript.
\bibitem  [Tei1]   {Tei1} Teissier B., \emph{The hunting of invariants in the geometry of discriminants.} In: P. Holm (ed.): Real and Complex Singularities, Oslo 1976, Northholland, 1978.
\bibitem  [Tei2]   {Tei2} Teissier B., \emph{R\'{e}solution simultan\'{e}e: II. R'{e}solution simultan\'{e}e et cycles \'{e}vanescents.} Lecture Notes in Math. \textbf{177}, 1980, 82-146.
\bibitem  [Wa]   {Wa} Wahl J., \emph{Equisingular deformations of plane algebroid curves.} Trans. Ams. Math. Sco. \textbf{193} (1974), 143-170.
\bibitem  [Za]   {Za} Zariski O., \emph{Contributions to the problem of equisingularity.} in "Questions on algebraic varieties", Ed. E. Marchionna, Cremonese, Roma, 1970, 261-343.
\end{thebibliography}
\end{document}